\documentclass{amsart}
  
\usepackage{amsthm}
\usepackage{amsmath}
\usepackage{amssymb}

\input xy  
\xyoption{all}

\author{Ian Shipman}
\date{\today}
\address{Dept. of Mathematics \\ University of Chicago \\ 5734 S. University Ave. \\ Chicago, IL 60637}
\email{ics@math.uchicago.edu}
\title{On representation schemes and Grassmanians of finite dimensional algebras and a construction of Lusztig}

\def\ZpI{\Z_{\geq 0}\{I\}}
\def\CI{\C\{I\}}
\def\bdim{\mb{dim}}
\def\Rep{\text{Rep}}
\def\Gr{\text{Gr}}
\def\v{\mb{v}}
\def\d{\mb{d}}
\def\dQ{\overline{Q}}

\def\A{\mathbb{A}}
\def\onto{\twoheadrightarrow}
\def\into{\hookrightarrow}
\def\tensor{\otimes}

\def\rvd{\Rep(A,V,D)}
\def\Rvd{\Rep^{st}(A,V,D)}
\def\grvd{\Gr^*_A(A_D,\mb{v})}
\def\sch{\mathfrak{sch}_\C}
\def\sets{\mathfrak{Set}}
\def\Gv{G_{\mb{v}}}
\def\iso{\stackrel{\simeq}{\longrightarrow}}
\def\Rvd{\Rep^{st}(A,V,D)}
\def\Ldv{\Lambda(D,V)^{st}}
\def\ldv{\mc{L}(D,V)}

\theoremstyle{plain}
\newtheorem{thm}{Theorem}[section]
\newtheorem{prop}[thm]{Proposition}
\newtheorem{cor}[thm]{Corollary}
\newtheorem*{mainthm}{Theorem \ref{maintheorem}}

\theoremstyle{definition}
\newtheorem{dfn}[thm]{Definition}
\newtheorem{rmk}[thm]{Remark}

\def\mc{\mathcal}
\def\frk{\mathfrak}

\def\mb{\mathbf}

\def\Z{\mathbb{Z}}
\def\C{\mathbb{C}}
\def\O{\mc{O}}
\def\onto{\twoheadrightarrow}
\def\into{\hookrightarrow}
\def\tensor{\otimes}
\DeclareMathOperator{\Hom}{Hom}
\DeclareMathOperator{\End}{End}
\DeclareMathOperator{\id}{id}

\def\rank{\text{rank}}
\def\dim{\text{dim}}

\begin{document}

\maketitle

\begin{abstract}
Let $I$ be a finite set and $\CI$ be the algebra of functions on $I$.  For a finite dimensional $\C$ algebra $A$ with $\CI \subset A$ we show that certain moduli spaces of finite dimsional modules are isomorphic to certain Grassmannian (quot-type) varieties.  There is a special case of interest in representation theory.  Lusztig defined two varieties related to a quiver and gave a bijection between their $\C$-points, \cite[Theorem 2.20]{L1}.  Savage and Tingley raised the question \cite[Remark 4.5]{ST} of whether these varieties are isomorphic as algebraic varieties.  This question has been open since Lusztig's original work.  It follows from the result of this note that the two varieties are indeed isomorphic.
\end{abstract}

\section{Introduction}
Fix a finite set $I$ and let $\CI$ be the algebra of functions on $I$.  It is a commutative semisimple algebra.  For each $i \in I$ there is an idempotent $e_i$, which satisfies $e_i(j) = \delta_{ij}$ for all $j \in I$.  Any right $\CI$ module $M$ is an $I$ graded vector space $M = \oplus_{i \in I}{ M_i }$ where $M_i = M\cdot e_i$.  If $M$ is finite dimensional we set $\bdim(M) = ( \dim(M_i) )_{i \in I} \in \ZpI$ and call elements of $\ZpI$ \emph{dimension vectors}.  For two dimension vectors $\mb{u},\v \in \ZpI$ we write $\mb{u} \leq \v$ if $u_i \leq v_i$ for all $i \in I$, with $\leq$ defined analogously.

Let $A$ be a finite dimensional $\CI$ over-ring, that is a finite dimensional $\C$ algebra with an embedding $\CI \subset A$.  We can view $A$ as an $I$ bi-graded algebra.  There are several natural moduli spaces associated to $A$.  First, after fixing a finite dimensional $\CI$ module $V$ we can form the representation scheme $\Rep(A,V)$ of isomorphism classes of right $A$ module structures on $V$ compatible with the $I$ grading.  From now on, we asssume that all $A$ modules are right $A$ modules and that an $A$ module structure on an $I$ graded vector space is compatible with the $I$ grading.  The scheme $\Rep(A,V)$ is a wild object and is properly studied in the setting of algebraic stacks.  To work with a more geometric moduli space, we can fix another finite dimensional $\CI$ module $D$ and form the ``$D$-framed'' representation scheme $\Rep(A,V,D)$.  This is the moduli spaces of right $A$ module structures on $V$ together with an arrow $D \to V$ that parameterizes a system of $A$-generators of $V$.  We construct $\Rep(A,V,D)$ using geometric invariant theory.  For a fixed $A$, the scheme $\Rep(A,V,D)$ only depends, up to isomorphism, on $\bdim(V)$ and $\bdim(D)$, so we can write either $\Rep(A,V,D)$ or $\Rep(A,\bdim(v),\bdim(D))$ to denote this scheme.

Next, suppose we fix a finite dimensional right $A$ module $M$ and a dimension vector $\v \in \ZpI$ with $\v \leq \bdim(M)$.  Then we have two more natural moduli spaces, the Grassmannians $\Gr_A(M,\v)$ and $\Gr_A^*(M,\v)$ of $\v$-dimensional submodules and quotients of $M$, respectively.

The purpose of this note is to compare certain of these constructions.  Let $D$ be a finite dimensional $\CI$ module.  Then $A_D = D \tensor_{\CI} A$ is a finite dimensional, projective right $A$ module.  Fix $\v \in \ZpI$ such that $\bdim(D) \leq \v \leq \bdim( A_D )$.  We prove
\begin{mainthm}
The two schemes $\Gr_A^*(A_D,\v)$ and $\Rep(A,\v,\bdim(D))$ are isomorphic. 
\end{mainthm}

There is an important special case.  If $Q$ is a quiver and $\C Q$ is the path algebra then we can form the preprojective algebra $\Pi_0$ and the quotients $A_N$ of this by powers of the augmentation ideal.  Then the moduli space of representations of $A_N$ is known as the nilpotent Lagrangian Nakajima quiver variety, and this result verifies that there is a quot-scheme description of it, as suggested by the work of Lusztig in \cite{L1}.  

\section{Main theorem}
In this section, we will recall the definitions and prove the main theorem.  For the following definition, we need only that $A$ is finitely generated over $\C$.  Let us start with the moduli space of framed representations.  Fix $\CI$ modules $D$ and $V$ of dimensions $\d,\v$ respectively.  We think of $\End_\C(V)$ as a $\CI$-bimodule in the obvious way, $e_i\phi e_j(v) = e_i \phi(e_j v)$.  If $R$ is any $\C$-algebra, then $\End_\C(V) \tensor_\C R$ inherits the bimodule structure.  Write $A^{op}$ for the algebra with the same underlying vector space as $A$ but with reversed multiplication.  Consider the functor 
\[ \Rep^{fr}(A,V)(R) = \Hom_{\CI-alg}(A^{op},\End_\C(V)\tensor_\C R) \]
from commutative $\C$-algebras to $\frk{S}ets$.  This functor associates to $R$ the set of $R$-linear right $A$ module structures on $V \tensor_\C R$.  This is co-represented by a finitely generated, commutative $\C$-algebra $\C[\Rep^{fr}(A,V)]$.  We denote by $\Rep^{fr}(A,V)$ the corresponding affine variety.  Let $\Gv = \prod_{i\in I} GL(v_i) = GL_{\CI}(V)$, the gauge group, and observe that $\Gv$ acts on the functor $\Rep^{fr}(A,V)$ through the conjugation action on $\End_\C(V)$, which preserves the bimodule structure.  Therefore it acts on the variety $\Rep^{fr}(A,V)$.  Let $\Rep^{fr}(A,V,D) = \Rep^{fr}(A,V) \times \Hom_{\CI}(D,V)$, equipped with the action of $\Gv$ on the factor $\Rep^{fr}(A,V)$ discussed above and the obvious action on $\Hom_{\CI}(D,V)$ by $\gamma(\phi) = \gamma \circ \phi$.
\begin{dfn}
The variety $\Rep(A,V,D)$ is defined to be the GIT quotient \cite{GIT} of $\Rep^{fr}(A,V,D)$ by $\Gv$ with respect to the inverse determinant character $\chi(\gamma) = \prod_i \det^{-1}(\gamma_i)$.
\end{dfn}
In order to form a GIT quotient, we first form the stable locus $\Rep^{st}(A,V,D)\subset \Rep^{fr}(A,V,D)$ with respect to $\chi$, which is an open set. Then we take the geometric quotient with respect to the action of $\Gv$.  This yeilds an arrow $\Rep^{st}(A,V,D) \onto \Rep(A,V,D)$.  

In what follows, it is necessary to identify $\Rep^{st}(A,V,D)$ and show that the action is free.  As preparation, observe that there is a tautological map $TA = T^{\tensor}_{\CI}{A} \onto A$, where $TA$ is the tensor algebra of $A$ viewed as a $\CI$ bimodule.  Now assume that $A$ is finite dimensional; then the tensor algebra is finitely generated over $\C$.  This map induces a $\Gv$ equivariant map of functors $\Rep^{fr}(A,V) \into \Rep^{fr}(TA,V)$ and thus a $\Gv$-equivariant morphism of algebraic varieties.  Clearly, 
\[ \Rep^{fr}(TA,V) \cong \Hom_{\CI-bimod}(A^{op},\End_\C(V)) = \bigoplus_{i,j \in I}{\Hom(e_i A e_j, \Hom(V_i, V_j))}.\]
Of course, there is also a closed embedding of varieties $\Rep^{fr}(A,V,D) \to \Rep^{fr}(TA,V,D)$.  Now, $\Rep^{st}(A,V,D)$ is just the inverse image of $\Rep^{st}(TA,V)$ under this embedding.  

Recall, that $\Rep^{st}(TA,V,D)$ is the locus of points $(x,p) \in \Rep^{fr}(TA,V,D)$ such that for any nonzero $t \in \C$, the orbit $\Gv \cdot ((x,p),t) \subset \Rep^{fr}(TA,V,D)\times \A^1$ is closed, where $\Gv$ acts on the factor $\A^1$ by the character $\chi$.  Any geometric point $x \in \Rep^{fr}(TA,V)$ defines a right $TA$ module structure on $V$, compatible with the $I$ grading, by $v \cdot a_1 \tensor \dotsm \tensor a_n = (x(a_n)\circ \dotsm \circ x(a_1))(v)$. As explained in \cite{N}, the point $((x,p),t)$ is stable if and only if the image under $p$ of $D$ generates $V$ as a $TA$ module with action map $x$.  Similarly, a geometric point $(x,p) \in \Rep^{fr}(A,V,D)$ is just a right $A$ module structure on $V$ together with a map $p:D \to V$ whose image generates $V$ over $A$.

Consider a geometric point $(x,p) \in \Rep^{st}(A,V,D)$, that is an $A$ module structure $x$ on $V$ and a map $p:D \to V$ whose image generates $V$.  An element $\gamma \in \Gv$ carries $(x,p)$ to $(\gamma x\gamma^{-1},\gamma p)$.  Suppose that $\gamma$ fixes $x$.  Then $\gamma$ is an $A$-module automorphism of $V$.  If $\gamma$ also fixes $p$ then $\gamma$ is the identity on a generating set for $V$ over $A$ and thus $\gamma = \id$.  Hence the action of $\Gv$ on $\Rep^{st}(A,V,D)$ is free.  We conclude that $\Rep^{st}(A,V,D) \onto \Rep(A,V,D)$ is a principal $\Gv$ bundle, or in other words a $\Gv$ torsor for the Zariski topology.

Next, fix a finite dimensional $A$ module $M$ and a dimension vector $\v$ with $\v \leq \bdim(M)$.  Then $\Gr_A^*(M,\v)$ is a subscheme of $\Gr^*(M,\v) := \prod_{i\in I}{\Gr^*(M_i,v_i)}$, where $\Gr^*(M_i,v_i)$ is the ordinary Grassmannian of $v_i$ dimensional quotient spaces of $M_i$.  Over $\Gr^*(M,\v)$ we have the $I$-graded Euler sequence
\begin{equation*}\label{euler_sequence}
 0 \to \mc{S} \to M \tensor_\C \O_{\Gr^*} \to \mc{Q} \to 0
\end{equation*}
and each element $a \in A$ defines an map $a:M \to M$ by the module structure.  Hence we can consider the arrow $x_a:\mc{S} \to \mc{Q}$ defined to be the composition
\[ \mc{S} \to M \tensor_\C \O_{\Gr^*} \stackrel{a\tensor \id}{\longrightarrow} M \tensor_\C \O_{\Gr^*} \to \mc{Q}. \]
\begin{dfn}[Grassmannian of quotient modules]
With respect to any local trivializations of $\mc{S}$ and $\mc{Q}$, each $x_a$ is represented by a matrix of functions.  We define $\Gr_A^*(M,\v)$ to be the subscheme of $\Gr^*(M,\v)$ whose ideal sheaf is generated, locally, by the entries of the matrices $x_a$ as $a$ varies through $A$.
\end{dfn}

Recall that $A_D = D \tensor_{\CI} A$.  In the remainder of this section, we will prove the following.
\begin{thm}\label{maintheorem}
The varieties $\Rep(A,V,D)$ and $\Gr^*_A(A_D,\v)$ are isomorphic.
\end{thm}

We work in the category $\sch$ of locally noetherian schemes over $\C$.  We will work in the Zariski topology below and this will be the context for the words ``local'' and ``locally''.  Also, we will use fpqc descent for morphisms and for quasicoherent sheaves, but the fpqc topology will not enter the picture in any other way.  By abuse of notation, $\rvd$ and $\grvd$ denote both the schemes defined above and the representable functors $\sch^{op} \to \sets$ that they define.  We will use the following folklore result and we include a proof for the convenience of the reader.

Let $G$ be an algebraic group over $\C$.  Suppose that $\pi_0:X \to Y$ is a principal $G$ bundle.  Given a locally noetherian scheme $T$, define a $T$-pair $(\pi,\phi)$ to be a principal $G$ bundle $\pi:P \to T$ and a $G$ equivariant map $\phi:P \to X$.  We say that two pairs $(\pi,\phi)$ and $(\pi',\phi')$ are equivalent if there are isomorphisms of principal $G$ bundles $\psi:P \to P'$ and $\psi_0:X \to X$ such that $\psi_0 \circ \phi = \phi' \circ \psi$.   Let $[X/G](T)$ be the set of $T$-pairs up to equivalence, i.e. the quotient stack of $X$ by $G$.  If $f:S \to T$ is a morphism then $[X/G](f):[X/G](T) \to [X/G](S)$ sends a pair $(\pi:P \to T,\phi)$ to $(f^*\pi:P \times_T S \to S,f^*\phi:P \times_T S \to X)$.  Given a morphism $f:T \to Y$, let $\pi_f:X \times_Y T \to T$ be the induced principal $G$ bundle and $\phi_f:X \times_Y T \to X$ the canonical $G$ equivariant morphism.
\begin{prop}\label{stacky}
The natural transformation $Y \to [X/G]$ defined by sending the morphism $f:T \to Y$ to the $T$-pair $(\pi_f,\phi_f)$ is an isomorphism of functors.
\end{prop}
\begin{proof}
We will construct an inverse to the map $Y(T) \to [X/G](T)$.  Let $(\pi,\phi)$ be a $T$-pair and set $g = \pi_0 \circ \phi:P \to Y$ be the composition.  Since the representable functor $Y$ is a sheaf in the fpqc topology (e.g. \cite[Exp. VIII, ch. 5]{SGA1}) and $P \to T$ is an fpqc cover, $g$ factors through a map $T \to Y$ if and only if $g$ equalizes the two projections $P \times_T P \rightrightarrows P$.  However since $P$ is a principal bundle $P \times_T P \cong G \times P$ and the two projections become the action $G \times P \to P$ and the second projection.  Now, $G$ equivariance means that we have a commutative diagram
\[
\xymatrix{
G \times P  \ar[r]^{\id \times \phi} \ar[d] & G \times X \ar[d] \\
P \ar[r]^{\phi} & X
}
\]
where the vertical arrows are the action morphisms.  Hence both of the arrows $G \times P \to Y$ factor through $G \times X \to Y$ and since $X \to Y$ is $G$ invariant, these are equal.  Hence $\phi$ can be factored $P \to T \to Y$ for some $f:T \to Y$ so $P$ is isomorphic to $f^*X$ and $\phi$ is isomorphic to the canonical equivariant morphism $f^*X \to X$.
\end{proof}

To get a better description of the functor $\rvd$, we will need a global notion of $A$ module.  Let $T$ be a locally noetherian scheme and let $\mc{F}$ be a locally free coherent sheaf.  A right $A$ module structure on $\mc{F}$ is a homomorphism $A^{op} \to \End(\mc{F})$ such that $\mc{F}_i := \mc{F} \cdot e_i$ is locally free for each $i \in I$.  We call such a structure a sheaf of $A$ modules or simply an $A$ module (over $T$).  The \emph{rank vector} $\bdim(\mc{F})$ is defined by $\bdim(\mc{F}) = (\rank(\mc{F}_i))_{i \in I} \in \Z_{\geq 0} I$.  If $\mc{F}$ and $\mc{F}'$ are two $A$ modules and $\alpha:\mc{F} \to \mc{F}'$ is a morphism of $\O_T$-modules, then we say $\alpha$ is a morphism of $A$ modules if for any $a \in A$, $\alpha \circ a = a \circ \alpha$ when we identify $a$ with its images in $\End(\mc{F})$ and $\End(\mc{F}')$ respectively.  Note that if $f:S \to T$ is any morphism then $f^*$ defines a functor from $A$ modules on $T$ to $A$ modules on $S$, thanks to the algebra maps $\End(\mc{F}) \to \End(f^*\mc{F})$.  This functor preserves rank vectors.

We will show that the following moduli-type functor is isomorphic to both $\rvd$ and $\grvd$.
\begin{dfn} Let $F$ be the presheaf of sets on the category of locally noetherian schemes defined by
\[ F(T) = \{\text{pairs } (\mc{V},p) \text{ on } T \}/\sim \]
where $\mc{V}$ is an $A$ module with $\bdim(\mc{V}) = \mb{v}$ and $p:D \tensor \mc{O}_T \to \mc{V}$ is an $I$ graded $\O_T$ module map such that in every geometric fiber of $\mc{V}$, the image of $p$ generates $\mc{V}$ as a $\Pi_0$ module.  The equivalence relation $\sim$ is defined as follows.  Two pairs $(\mc{V},p), (\mc{V}',p')$ are equivalent if and only if there exists an isomorphism $\gamma:\mc{V} \to \mc{V}'$ of $A$ modules such that $p'_i = \gamma_i p_i$.  If $f:S \to T$ is a morphism then $F(f)[(\mc{V},p)] = [(f^*\mc{V}, f^*p)]$, where we use the canonical identification $D \tensor_\C \O_S \cong f^*( D \tensor_\C \O_T)$
\end{dfn}

First, we construct a natural transformation $F \Rightarrow \rvd$.  Let $T$ be a locally noetherian scheme and let $(\mc{V},p)$ represent a class in $F(T)$.  Let $P_\mc{V}$ be the $I$-graded frame bundle of $\mc{V}$.  Explicitly, if $P_{\mc{V}_i}$ is the ordinary frame bundle of $\mc{V}_i$ then $P_\mc{V} = \times_I P_{\mc{V}_i}$, fiber product over $T$.  By construction, $\pi:P_\mc{V} \to T$ is a $\Gv$ torsor for the Zariski topology which has the property that $\pi^*\mc{V}$ is a $\Gv$-equivariantly trivial $I$-graded vector bundle.  Choose a $\Gv$-equivariant trivialization $V \tensor_\C \O_{P_\mc{V}} \iso \pi^*\mc{V}$.  This endows $V \tensor_\C \O_{P_\mc{V}}$ with the structure of an $A$ module.  Composing the pullback of $p$ with the trivialization gives a map $p:D \tensor_C \O_{P_\mc{V}} \to V \tensor_\C \O_{P_\mc{V}}$ which generates $V \tensor_\C \O_{P_\mc{V}}$ in every geometric fiber.  The $A$ module structure and map $p$ are just matrix valued functions on $P_\mc{V}$ satisfying certain identities.  They define a morphism into $\Rvd$, which is $\Gv$ equivariant.  By Proposition \ref{stacky}, the pair $(P_\mc{V} \to T, P_\mc{V} \to \Rvd)$ defines a map $T \to \rvd$.  It is straightforward to check that this construction does not depend on any of the choices that we made and is functorial.  

Next, we construct the inverse natural transformation $\rvd \Rightarrow F$.  Let $f:T \to \rvd$ be a morphism.  There is an important canonical element $\gamma_0$ of $F(\rvd)$.  We construct this in the following way.  There is a canonical $A$ module $V \tensor \O_{\Rvd}$ on $\Rvd$.  It has a $\Gv$-equivariant structure which is compatible with the $A$ module structure, meaning that each $a \in A$ acts on $V \tensor_\C \O_{\Rvd}$ by a morphism of $\Gv$-equivariant sheaves.  The $\Gv$-equivariant pair $(V\tensor_\C \O_{\Rvd},\tilde{p}_0:D \tensor_\C \O_{\Rvd} \to V\tensor_\C \O_{\Rvd})$ descends (e.g. fpqc descent as in \cite[Exp. VIII, ch. 1]{SGA1}) to a pair $(\mc{V}^0,p_0)$ on $\rvd$.  We call the class of this pair $\gamma_0$.  We send $f$ to $f^*\gamma_0 \in F(T)$.  It is straightforward to check that these two natural transformations are mutually inverse.

Now, we construct a natural transformation $\grvd \Rightarrow F$.  By construction of $\grvd$, a map $T \to \grvd$ is the same thing as an $I$-graded quotient bundle $\mc{Q}$ fitting into an Euler-type exact sequence
\[ 0 \to \mc{S} \to A_D \tensor_\C \O_T \to \mc{Q} \to 0\]
having rank vector $\mb{v}$ and such that certain maps $\mc{S}_i \to \mc{Q}_j$ vanish.  The bundle $\mc{S}$ is defined by the surjection $A_D \tensor_\C \O_T \to \mc{Q}$.  Since these maps vanish, $\mc{Q}$ and $\mc{S}$ are $\Pi_0$ sub- and quotient modules of $A_D \tensor_\C \O_T$.  We form $p: D \tensor_\C \O_T \to \mc{Q}$ by composing the natural inclusion $D \tensor_\C \O_T \into A_D \tensor_\C \O_T$ with the projection to $\mc{Q}$.  Since $A_D \tensor_\C \O_T$ is generated by the image of $p$ in each fiber, so is $\mc{Q}$.  Hence $(\mc{Q},p)$ determines an element in $F(T)$.  Once again it is straightforward to check that this construction is functorial.

Finally, we construct the inverse natural transformation $F \Rightarrow \grvd$.  Given a pair $(\mc{V},p)$ in $F(T)$ we can extend $p$ to a map $A_D \tensor_\C \O_T \to \mc{V}$.  Since $T$ is locally noetherian and this map is fiberwise surjective, it is surjective.  Thus we get an exact sequence
\[
0 \to S \to A_D \tensor_\C \O_T \to \mc{V} \to 0
\]
where $S$ is simply the kernel.  Now since the quotient map is graded of constant rank, $S$ is an $I$-graded vector bundle and since $\mc{V}$ is a quotient $A_D$ module, those maps $S_i \to \mc{V}_j$ vanish.  Therefore the quotient $A_D \tensor_\C \O_T \onto \mc{V}$ determines a map $T \to \grvd$.  This natural transformation is inverse to the one constructed in the previous paragraph.

In summary, there are isomorphisms $\rvd \Longleftrightarrow F \Longleftrightarrow \grvd$.  Yoneda's lemma implies that the composition is induced by mutually inverse algebraic maps between $\grvd$ and $\rvd$.  This proves our result.

\begin{rmk}\label{C_points}
The geometric points of $\rvd$ are isomorphism classes of pairs of an $A$ module structure on $V$ and a choice $p:D \to V$ of a generating set.  The geometric points of $\grvd$ are $A$ module structures on $V$ together with surjective $A$ module maps $A_D \onto V$.  On the level of geometric points the identification above becomes the following.  Given a geometric point of $\rvd$ we extend $D \to V$ to a map $A_D \to V$.  Because $D$ generates, this map is surjective.  Given a geometric point of $\grvd$ we keep the $A$ module structure on $V$ and define $p:D \to V$ to be the composite, $D \into A_D \onto V$.
\end{rmk}

\section{Quiver varieties}
The following notions are essentially standard, and are explained in \cite{N} under a different system of notation.  A quiver is a finite directed graph described by a pair $Q = (I,E)$ where $I$ is the set of vertices and $E$ is the set of edges.  For an edge $a \in E$, $a_\tau$ and $a_\sigma$ will denote the target and source vertices, respectively.
\[
a_\sigma \stackrel{a}{\longrightarrow} a_\tau
\]  
Fix a quiver $Q$ for the rest of the section and assume that $Q$ has no edge-loops, edges with the same source and target vertex.  Let $\C E$ be the vector space generated by the edges of $Q$.  There is a natural $\CI$ bimodule structure on $\C E$ determined by
\[ 
e_i a e_j =
\begin{cases}
a & a_\sigma = i, a_\tau = j, \\
0 & \text{otherwise}
\end{cases}
\]
Now, one defines the \emph{path algebra} of $Q$ to be $\C Q := T^{\tensor}_{\CI}{\C E}$, the tensor algebra over $\CI$ of the bimodule $\C E$.  The set of paths through $Q$ can be identified with a basis of $\C Q$ and under this identification, multiplication becomes concatenation of paths or zero if the paths cannot be concatenated.  The subalgebra $\CI \subset \C Q$ corresponds to the paths of length zero.  See \cite{L1,N,G} for further discussion.

A representation of $Q$ is a right $\C Q$ module.  If $M$ is a finite dimensional right $\C Q$ module then the dimension vector $\bdim(M) = (\dim_\C M\cdot e_i )_{i \in I} \in \Z_{\geq 0} I$ is just the dimension vector of $M$ viewed as a $\CI$ module.  If $V$ is a $\CI$ module, then $\C Q$ module structures extending the $\CI$ module structure are determined by a collection of linear maps $\{x_a:V_{a_\sigma} \to V_{a_\tau}\}_{a \in E}$.

The path algebra $\C Q$ is naturally graded by path length.  Let $\C Q_n$ be the subspace spanned by paths of length $n$ and let $\C Q_{\geq n} = \bigoplus_{m \geq n} \C Q_m$.  Note that $\C Q_{\geq n}$ is a two sided ideal.  Following Lusztig, we consider nilpotent representations.  Say that a representation $M$ of $Q$ is \emph{nilpotent} if $M \cdot \C Q_{\geq n} = 0$ for some $n$.  Let $V$ be a finite dimensional representation of $Q$ with $N = \dim_\C V$.  It is a consequence of \cite[1.8]{L2} that if $V$ is nilpotent then $V \cdot \C Q_{\geq N} = 0$.  Hence the module structure on $V$ factors through $\C Q / \C Q_{\geq N}$.

Let $V$ be a finite dimensional $\CI$ module with $\mb{v} = \bdim(V)$.  The \emph{representation space} of $Q$ on $V$ is
\[ \Rep(Q,V) = \bigoplus_{a \in E}{ \Hom_\C(V_{a_\sigma},V_{a_\tau})}. \]
In the notation of the previous section, $\Rep(Q,V) = \Rep^{fr}(\C Q, V)$.  Let $\Rep^{nil}(Q,V)$ be the space of nilpotent representations of $Q$.  Explicitly, $\Rep^{nil}(Q,V)$ is the subvariety of $\Rep(Q, V)$ defined by the matrix equations $x_{a^1} \dotsm x_{a^N} = 0$ over all paths $a^1 \dotsm a^N$ of length $N$, by the remark in the previous paragraph.  Since any nilpotent right $\C Q$ module structure on $V$ factors through $\C Q / (\C Q )_{\geq N}$ we have $\Rep^{nil}(Q,V) = \Rep^{fr}(\C Q / (\C Q )_{\geq N}, V)$.

Given a quiver $Q = (I,E)$, let $Q^{op} = (I,E^{op})$ be the \emph{opposite quiver} where $E^{op} = \{a^* | a \in E \}$ and the edges are attached by $a^*_\sigma = a_\tau$ and $a^*_\tau = a_\sigma$.  
\[  a_\sigma \stackrel{a^*}{\longleftarrow} a_\tau \]
Finally, the \emph{double quiver} is defined to be $\dQ = (I, E \sqcup E^{op})$.  In some presentations, e.g. \cite{L1}, one begins with an unoriented graph and the double quiver appears once one fixes a choice of orientation this underlying graph.  The path algebra of the double quiver has a distinguished element $\theta = \sum_{a \in E}{[a,a^*]}$, where $[a,a^*] = aa^* - a^*a$.  The quotient $\Pi_0 = \Pi_0(\dQ) = \C \dQ / (\theta)$ is an associative algebra called the \emph{preprojective algebra}.  

We have $\Rep^{fr}(\Pi_0,V) \subset \Rep(Q,V)$, embedded as the subscheme defined by the matrix equation $\sum_{a \in E}{[x_a,x_{a^*}]} = 0$.  The $\C$ points of this subscheme define $\Pi_0$ module structures on $V$ extending the $\CI$ module structure.  There is also the corresponding subscheme of nilpotent representations of $\Pi_0$, $\Rep^{nil}(\Pi_0,V) \subset \Rep(\Pi_0,V)$.  All of these schemes are $G_\mb{v}$ invariant. 

\begin{dfn}[Lagrangian Nakajima quiver variety]
Let $\Lambda(D,V) = \Rep^{nil}(\Pi_0,V)\times \Hom_{\CI}(D,V)$.  The action of $G_\mb{v}$ on $\Rep^{nil}(\Pi_0,V)$ extends to $\Lambda(D,V)$ by
\[ \gamma (x_a, p_i) = ( \gamma_{a_\tau} x_a \gamma_{a_\sigma}^{-1}, \gamma_i p_i) \]
where $\gamma = (\gamma_i)_{i \in I}$ and $\{p_i:D_i \to V_i\}$ are the components of a map $p \in \Hom_{\CI}(D,V)$.  There is an open set $\Ldv \subset \Lambda(D,V)$ of stable points.  The open set is defined by the GIT stability condition \cite{GIT} associated to the character $\chi(g_i) = \prod_i \det(g_i)^{-1}$.  As explained in \cite{N}, if we think of a point $(x_a,p_i)$ as defining a $\Pi_0$ module structure on $V$ and a map $D \to V$ then this point is in $\Ldv$ if and only if the image of $p$ generates $V$ as a $\Pi_0$ module.  By \cite[\S 1]{GIT}, the geometric quotient scheme $\ldv := \Ldv/G_\mb{v}$ exists and since $G_\mb{v}$ acts freely on $\Ldv$, the natural morphism $\Lambda(D,V) \to \ldv$ is a principal $G_\mb{v}$ bundle.  Note that the Lagrangian Nakajima quiver variety $\ldv$ only depends, up to isomorphism, on $\mb{d}$ and $\mb{v}$. In the formalism of the previous section, $\ldv = \Rep(\Pi_0/(\Pi_0)_{\geq N},V,D)$.
\end{dfn}

Let $\mb{P}_D = D \tensor_{\CI} \Pi_0$ and consider $\Gr^*_{\Pi_0}(\mb{P}_D,V)^{nil}$, the Grassmannian of nilpotent quotients of $\mb{P}_D$ of dimension $\bdim(V)$.  We can also think of $\Gr^*_{\Pi_0}(\mb{P}_D,V)^{nil}$ as the moduli space of nilpotent right $\Pi_0$ module structures on $V$ and surjections $\mb{P}_D \onto V$ of right $\Pi_0$ modules.  Let $\mb{P}^N_D = D \tensor_{\CI} \Pi_0/(\Pi_0)_{\geq N} = \mb{P}_D/ \mb{P}_D \cdot (\Pi_0)_{\geq N}$.  Since any nilpotent right $\Pi_0$ module structure on $V$ factors through $\Pi_0/(\Pi_0)_N$ we have $\Gr^*_{\Pi_0}(\mb{P}_D,V)^{nil} \cong \Gr^*_{\Pi_0/(\Pi_0)_{\geq N}}(\mb{P}^N_D,V)$.  We now have the following consequence of \ref{maintheorem}.
\begin{cor}
 The schemes $\mc{L}(V,D)$ and $\Gr^*_{\Pi_0}(\mb{P}_D, V)^{nil}$ are isomorphic.
\end{cor}

It is clear from remark \ref{C_points} that the induced bijection on $\C$ points is just the one Lusztig gave in \cite{L1}.

\subsubsection*{Acknowledgements}  
I would like to thank Victor Ginzburg for bringing this problem to my attention.  Also, I wish to thank Tom Church for pointing out that the frame bundle can be used to trivialize a vector bundle.  Also, special thanks to Peter Tingley for carefully reading a draft version of this article and offering useful comments and discussion.

\bibliographystyle{plain}
\bibliography{quiver_varieties}

\end{document}